\documentclass[12pt]{article}

\title{The Jammed Phase of the
Biham-Middleton-Levine Traffic Model}
\author{Omer Angel\thanks{Funded
in part by the Pacific Institute of Mathematical Sciences.},
Alexander E. Holroyd\thanks{Funded in part by an NSERC (Canada)
grant.} and James B. Martin}
\date{March 30, 2005 (revised June
29, 2005)}

\usepackage{amsmath,amsthm,amssymb,graphics,color}

\newtheorem{thm}{Theorem}
\newtheorem{lemma}[thm]{Lemma}
\newtheorem{proposition}[thm]{Proposition}

\newtheorem{claim}[thm]{Claim}

\renewenvironment{proof}
{\begin{trivlist} \item \noindent{\sc Proof. }}
{\qed \end{trivlist}}
\newenvironment{proofof}[1]
{\begin{trivlist} \item \noindent{\sc Proof of\ #1. }}
{\qed \end{trivlist}}

\newcommand{\PP}{\mathbb{P}}
\newcommand{\ZZ}{\mathbb{Z}}

\newcommand{\ba}{\mathbf{a}}
\newcommand{\bb}{\mathbf{b}}
\newcommand{\bu}{\mathbf{u}}
\newcommand{\bv}{\mathbf{v}}
\newcommand{\bx}{\mathbf{x}}
\newcommand{\by}{\mathbf{y}}
\newcommand{\bz}{\mathbf{z}}

\newcommand\rta{\ensuremath{\boldsymbol{\rightarrow}}}
\newcommand\upa{\ensuremath{\boldsymbol{\uparrow}}}

\begin{document}
\maketitle
\renewcommand{\thefootnote}{}

\begin{abstract}
Initially a car is placed with probability $p$ at each site of the
two-dimensional integer lattice. Each car is equally likely to be
East-facing or North-facing, and different sites receive
independent assignments. At odd time steps, each North-facing car
moves one unit North if there is a vacant site for it to move
into. At even time steps, East-facing cars move East in the same
way. We prove that when $p$ is sufficiently close to $1$ traffic
is jammed, in the sense that no car moves infinitely many times.
The result extends to several variant settings, including a model
with cars moving at random times, and higher dimensions.
\footnote{\hspace{-2em}This work was done while all three authors were
hosted by MSRI in Berkeley.}
\footnote{\hspace{-2em}{\bf Key words:} traffic, phase transition,
percolation, cellular automata}
\footnote{\hspace{-2em}{\bf 2000 Mathematics Subject Classifications:}
Primary 60K35; Secondary 82B43}
\end{abstract}

\section{Introduction}

\begin{figure}
\centering \vspace{-0.3in}
\resizebox{2.2in}{!}{\includegraphics[0,0][200,200]{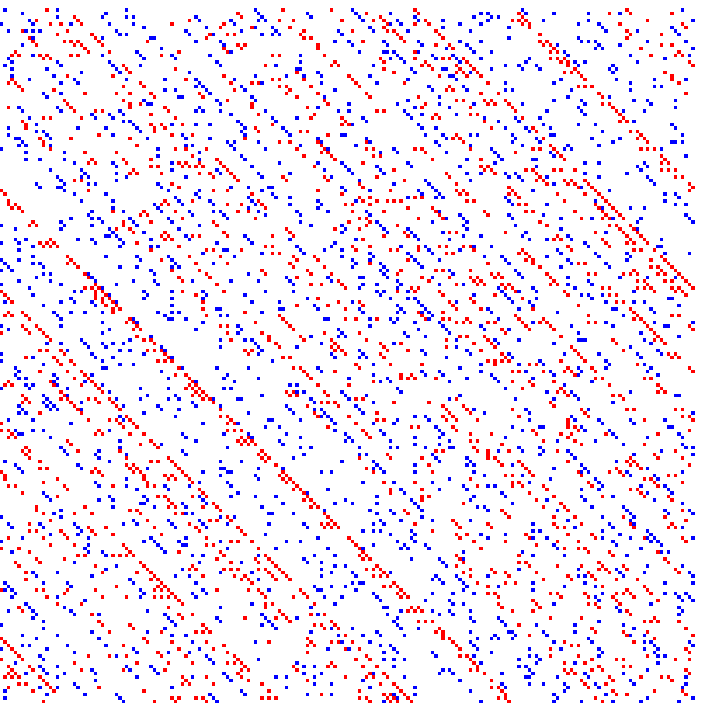}}
\resizebox{2.2in}{!}{\includegraphics[0,0][200,200]{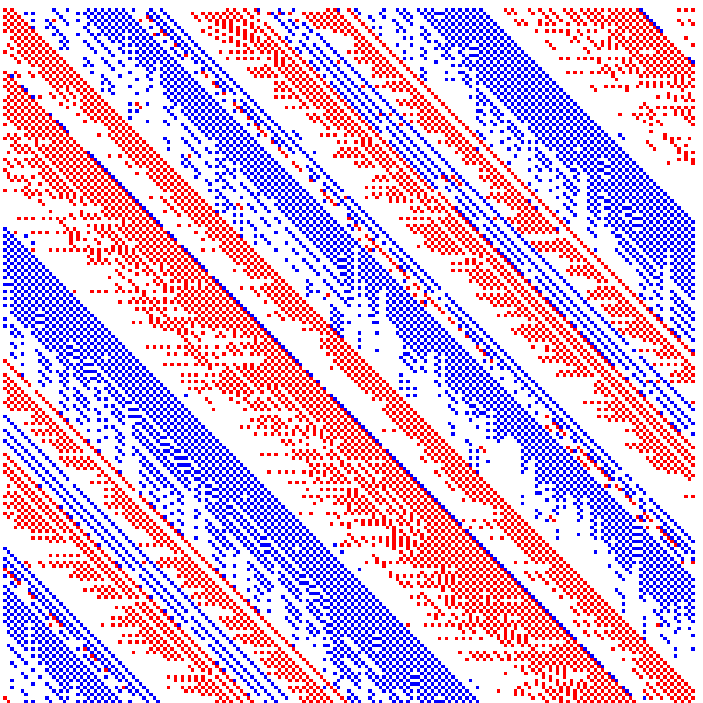}}\\
\vspace{-0.05in}
$p=0.1$ (free-flowing) \qquad\qquad $p=0.3$ (free-flowing)\\
\vspace{0.05in}
\resizebox{2.2in}{!}{\includegraphics[0,0][200,200]{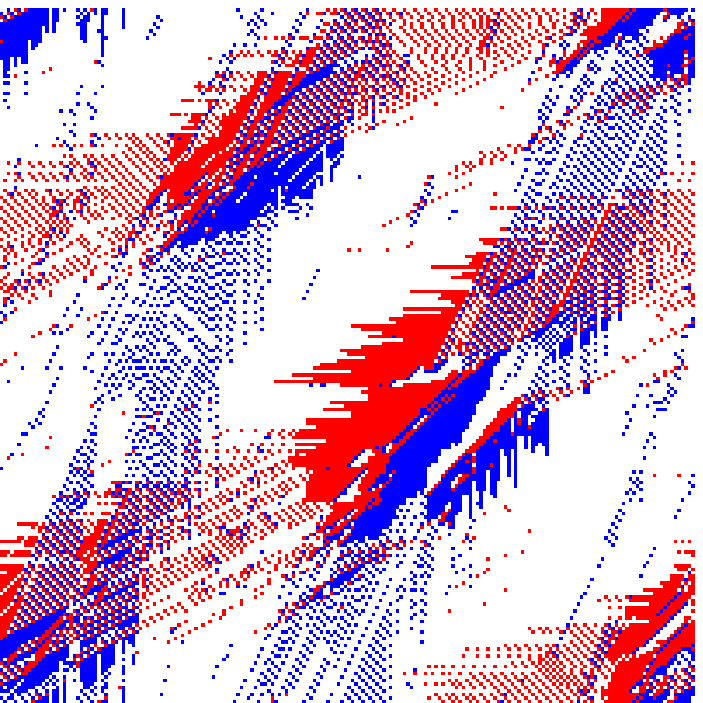}}
\resizebox{2.2in}{!}{\includegraphics[0,0][200,200]{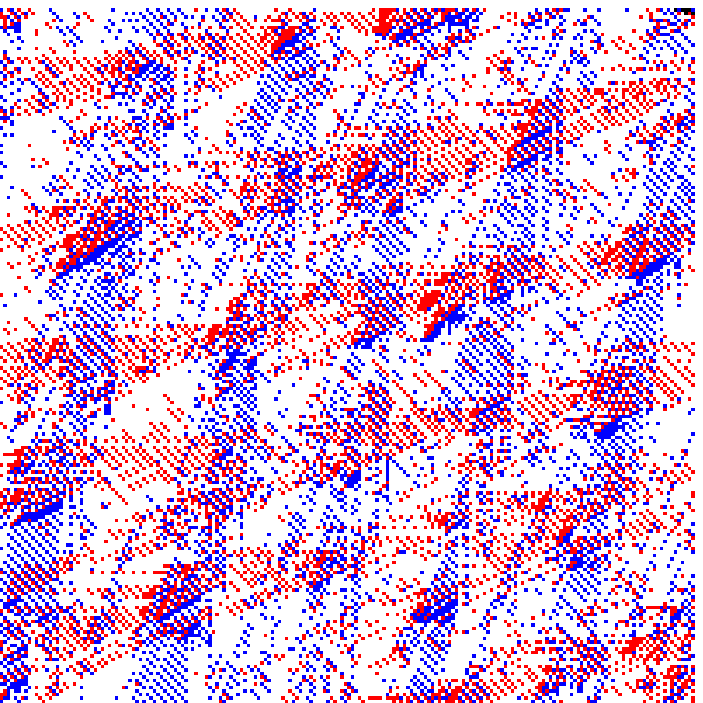}}\\
\vspace{-0.05in}
$p=0.32$ (intermediate) \qquad\qquad $p=0.32$ (intermediate)\\
\vspace{0.05in}
\resizebox{2.2in}{!}{\includegraphics[0,0][200,200]{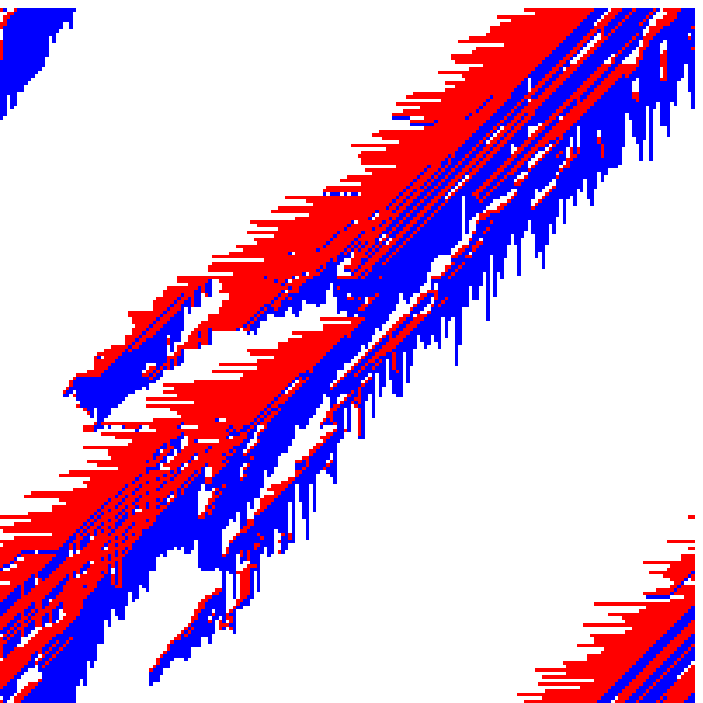}}
\resizebox{2.2in}{!}{\includegraphics[0,0][200,200]{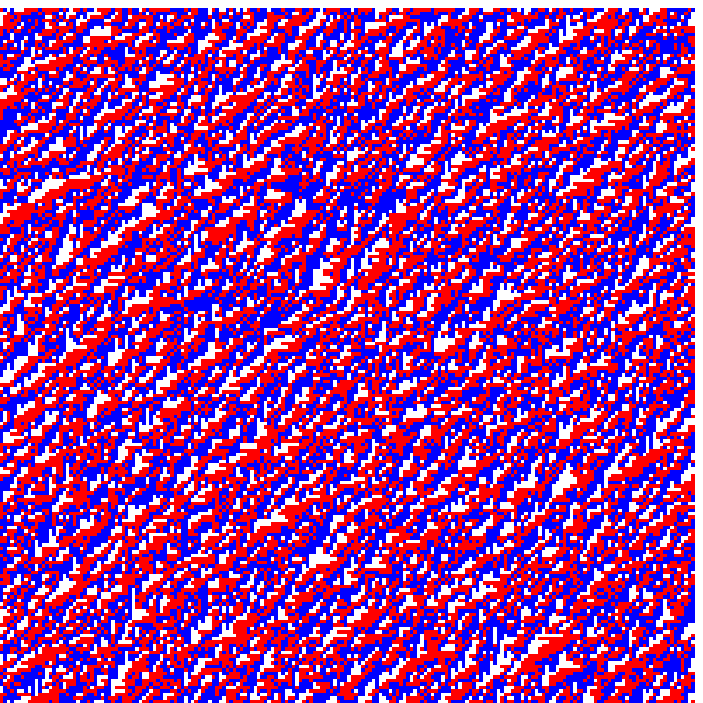}}\\
\vspace{-0.05in}
$p=0.34$ (jammed) \qquad\qquad $p=0.8$ (jammed)\\
\vspace{0.05in}
\caption{Examples of the model after 20,000 steps on a 200-by-200
torus. East-facing and North-facing cars are shown in red and blue
respectively.} \label{fig:pics}
\end{figure}

The following simple model for traffic congestion was introduced in
\cite{bml}.  Let $\ZZ^2=\{\bz=(z_1,z_2):z_1,z_2\in\ZZ\}$ be the
two-dimensional integer lattice.  At each time step $t=0,1,\ldots$,
each site of $\ZZ^2$ contains either an East car ($\rta$),
a North car ($\upa$) or an empty space (0).  Let $p\in[0,1]$.  The
{\em initial} configuration is given by a random element $\sigma$ of
$\{0,\rta,\upa\}^{\ZZ^2}$ under a probability measure $\PP_p$ in
which
\[
  \PP_p(\sigma(\bz)=\;\rta) = \PP_p(\sigma(\bz)=\;\upa) = p/2
 \quad\text{and}\quad \PP_p\big(\sigma(\bz)=0\big) = 1-p
\]
for each site $\bz\in\ZZ^2$, and the initial states of different sites are
independent.

The configuration evolves in discrete time according to the following
deterministic dynamics. On each odd time step, every $\upa$ which
currently has a $0$ immediately to its North (i.e. in direction (0,1))
moves into this space.  On each even time step, each $\rta$ which
currently has a $0$ immediately to its East (i.e. in direction (1,0))
moves into this space.  The configuration remains otherwise unchanged.

\begin{thm}
\label{maintheorem} There exists $p_1<1$ such that for all $p\geq
p_1$, almost surely no car moves infinitely often and the state of
each site is eventually constant.
\end{thm}

The above result goes part way towards establishing the
following natural conjecture.

\paragraph{Conjecture.}  There exists $p_c\in(0,1)$ such that for
$p>p_c$ almost surely no car moves infinitely often, while for
$p<p_c$ almost surely all cars move infinitely often.
\vspace{\parsep}

The model may be defined on a finite torus (i.e. a rectangle with
periodic boundary conditions) in a natural way, and our proof can
be adapted to this case.

\begin{thm} \label{torus}
  Consider the model on an $m$ by $n$ torus. There exists $p_1<1$ such that
  for any $p\geq p_1$ and any sequence of tori such that $m,n\to\infty$ and
  $m/n$ converges to a limit in $(0,\infty)$, asymptotically almost surely
  no car moves infinitely often and the configuration is eventually
  constant.
\end{thm}

The present article represents the first rigorous progress on the
model, which has previously been studied extensively via
simulation and partly non-rigorous methods.  Such studies have
suggested that $p_c\approx 0.35$, and furthermore that for $p$
sufficiently small (and perhaps even for {\em all} $p<p_c$),
 all cars move with asymptotic speed equal to the maximum
possible ``free flowing'' speed of $1/2$.  The latter striking
phenomenon was observed experimentally in \cite{bml}, and has been
conjectured for the infinite lattice by Ehud Friedgut (personal
communication).  Recent results in \cite{souza} suggest the
existence of further intermediate phases (involving speeds
strictly in $(0,1/2)$) for the model on finite tori.  The model
appears to exhibit remarkable self-organizing behaviour. The
problem of rigorously analyzing the model was given as an
``unsolved puzzle'' in \cite{winkler}. References to earlier work
may be found in \cite{souza}.  Some simulations are illustrated in
Figure~\ref{fig:pics}.

Here is an overview of our proof of Theorem \ref{maintheorem}. First
consider the (trivial) case $p=1$. Any given car is blocked by another car
immediately in front of it, this car in turn is blocked by a further car,
and so on. Thus the original car can never move because there is an
infinite chain of cars blocking it. This argument does not extend to $p<1$
because such a chain will always be broken by an empty space. Returning to
the case $p=1$, we therefore consider an additional local configuration
which can cause a car to be blocked in a different way, and which gives
rise to additional types of blocking paths. This local configuration occurs
with some positive intensity throughout space, therefore (for $p=1$) we
obtain an extensive network of blocking paths, any of which block a given
car. Taking $p<1$ is the same as removing a proportion of cars from a $p=1$
configuration. If the proportion of cars removed is sufficiently small, it
is likely that some of the blocking paths will survive, in which case the
original car will be blocked even when $p<1$. This argument is formalized
via a comparison with super-critical oriented percolation on a renormalized
lattice.

In principle, our arguments give an explicit bound for $p_1$ in
Theorem \ref{maintheorem}.  We do not attempt to compute this bound,
since it would be very close to 1, and nowhere near the supposed value
of $p_c$.

Our proof extends to yield analogous results in a number of variant
settings. These include: a model in which cars move at random Poisson
times rather than at alternate discrete time steps, initial conditions
with different probabilities of East and North cars, and higher
dimensional generalizations.  We discuss these variants, and Theorem
\ref{torus} concerning the torus, at the end of the article.

\section{Proof of Main Result}

\begin{figure}
\centering
\rotatebox{-90}{\reflectbox{\resizebox{!}{1.5in}{\includegraphics{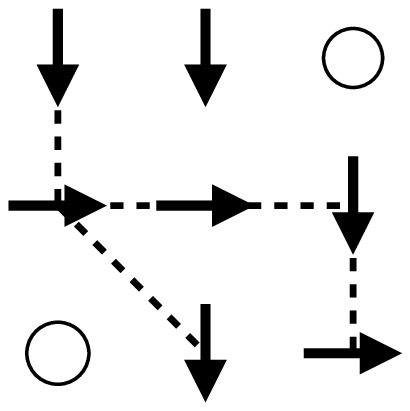}}}}
\caption{There are blocking paths from (0,0) (bottom left) to
(2,2) and from (0,0) to (1,2).  The latter uses a step of type
(iii).} \label{fig:block}
\end{figure}

A finite or infinite sequence of sites $\bz^0, \bz^1, \bz^2, \dots
[,\bz^n]$ is called a {\bf blocking path} if, for each $m\geq 0$,
one of the following holds:
%
\begin{list}{}{\setlength{\leftmargin}{.55in}}
\item[(i)] $\sigma(\bz^m)=\;\rta$ and $\bz^{m+1}=\bz^m+(1,0)$;
\item[(ii)]  $\sigma(\bz^m)=\;\upa$ and $\bz^{m+1}=\bz^m+(0,1)$;
\item[(iii)] $\sigma(\bz^m)=\sigma(\bz^m+(1,0))=\;\rta$, \quad
$\sigma(\bz^m+(1,-1))=\;\upa$, \\ and
$\bz^{m+1}=\bz^m+(1,1)$;
\item[or (iv)] $\sigma(\bz^m)=\sigma(\bz^m+(0,1))=\;\upa$, \quad
$\sigma(\bz^m+(-1,1))=\;\rta$, \\ and
$\bz^{m+1}=\bz^m+(1,1)$.
\end{list}
See Figure \ref{fig:block} for an illustration.  Note that if
$\bz^0, \ldots, \bz^n$ and $\bz^n,\bz^{n+1},\ldots$ are blocking
paths then so is $\bz^0, \ldots, \bz^n,\bz^{n+1},\ldots$.  Cases (i) and
(ii) correspond to the na\"{\i}ve chains
of cars mentioned in the introduction.  Cases (iii) and (iv) will provide the key to our argument by allowing for additional types of blocking path.

\begin{lemma}
\label{blocking}
No car on an infinite blocking path ever moves.
\end{lemma}

\begin{proof}
We claim that the car at $\bz^m$ can only move strictly after that
at $\bz^{m+1}$ has moved.  This implies the result, by induction
on the time step.  The claim is immediate in cases (i) and (ii)
above.  In case (iii), we note that the car at $\bz^m$ can only
move after that at $\bz^m+(1,0)$.  If the latter car ever moves
then it does so at an even step, and it is replaced immediately at
the next step by the car initially at $\bz^m+(1,-1)$.  But this
car now cannot move again until after that at $\bz^{m+1}$.  An
analogous argument applies in case (iv).
\end{proof}

\begin{figure}
\centering \resizebox{!}{2.5in}{\includegraphics{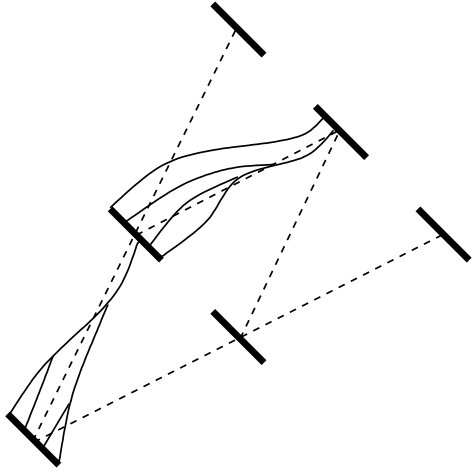}}
\caption{Part of the renormalized lattice. Renormalized sites are indicated
  by bold lines, renormalized edges by dashed lines and blocking paths by
  curved lines. Here $((0,0),(0,1))$ and $((0,1),(1,1))$ are good edges.}
\label{fig:renorm}
\end{figure}

We introduce a renormalized lattice with the structure of $\ZZ^2$.
Let $M,k$ be integers (to be fixed later) satisfying $M>2k>0$.
 Each site in the
renormalized lattice consists of $2k+1$ sites on a diagonal. Denote by
$D_k$ the set $\{(s,-s) : |s|\le k\}$. For each site
$\bu=(u_1,u_2)\in\ZZ^2$ we define the renormalized site
\[
  V_\bu = u_1(10M,9M) + u_2(9M,10M) + D_k.
\]
A renormalized {\bf edge} is an ordered pair $(\bu,\bv)$ where
$\bv-\bu$ equals $(1,0)$ or $(0,1)$.  We say that the edge $(\bu,\bv)$
is {\bf good} if
\[
\text{from every $\bx\in V_{\bu}$ there is a blocking path to some
} \by\in V_{\bv}.
\]
(Recall that blocking paths, and therefore good edges, are defined in terms
of the initial configuration $\sigma$). See Figure \ref{fig:renorm} for an
illustration.

\begin{lemma}\label{mdeplemma}
Suppose $M>2k>0$.  The process of good edges is {\em
30-dependent}. (That is, if $A,B$ are sets of edges at graph-theoretic
distance at least 30 from each other in the renormalized lattice, then
the states of the edges in $A$ are independent of those in $B$).
\end{lemma}

\begin{proof}
From the definitions of blocking paths and good edges, the event
that the edge $(\bu,\bv)$ is good depends only on the initial
states $\sigma(\bx)$ of sites $\bx$ in a certain box containing
$V_{\bu}$ and $V_{\bv}$.  Since $M>2k$, such boxes are disjoint
for edges at graph-theoretic distance at least $30$ from each
other.
\end{proof}

\begin{figure}
\centering \resizebox{!}{3in}{\includegraphics{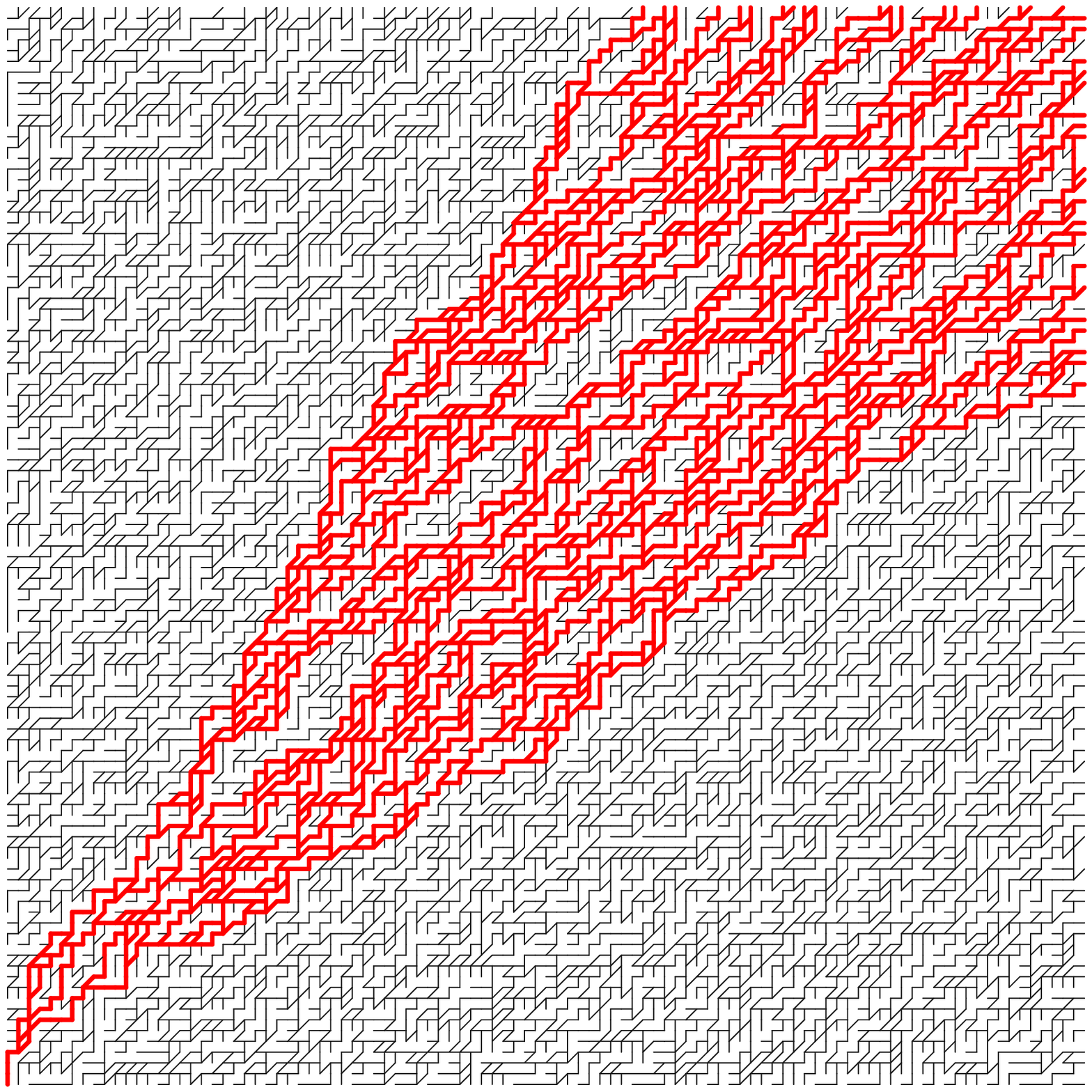}}
\caption{Blocking paths for a random initial configuration with $p=1$.
Blocking paths from the origin are highlighted.}
\label{fig:cone}
\end{figure}

In order to prove Theorem \ref{maintheorem} we will show that the
probability that an edge is good is close to 1. We will do this
first for the case $p=1$. Figure~\ref{fig:cone} illustrates all
blocking paths starting at the origin for a random initial
configuration with $p=1$. A key step is the following lemma which
states that such paths are likely to come close to any site in a
certain cone.

\begin{lemma} \label{target}
  Consider the case $p=1$. Let $E(\by,k)$ be the event that there is a
  blocking path from $(0,0)$ to $\by+(s,-s)$ for some $s\in[-k,k]$. There
  exists $c>0$ such that for any site $\by = (y_1,y_2) \in\ZZ^2$ satisfying
  $y_1,y_2>0$ and $y_1/y_2 \in [8/9,9/8]$ we have
  \[
  \PP_1\big(E(\by,k)\big) > 1-e^{-ck}.
  \]
\end{lemma}

The proof of Lemma~\ref{target} is deferred to the end of the section.

\begin{proposition}\label{goodboxprop2}
Consider the case $p=1$. For any $\beta<1$, there exist $M$ and $k$ with
$M>2k$ such that for every renormalized edge $(\bu,\bv)$,
\[
\PP_1(\text{\rm edge $(\bu,\bv)$ is good})\geq \beta.
\]
\end{proposition}

\begin{proof}
Take $k$ large enough that $(2k+1)e^{-ck} < 1-\beta$, and then take $M>2k$
large enough that
\[
  \frac{10M+k}{9M-k}\leq\frac{9}{8}.
\]
By Lemma~\ref{target} and translation invariance we obtain
\begin{align*}
  \PP_1(\text{edge $(\bu,\bv)$ is not good})  & \leq \sum_{\bx\in V_\bu}
  \PP_1(\nexists\text{ a blocking path from $\bx$ to $V_\bv$}) \\
   & \leq (2k+1)e^{-ck} < 1-\beta.
\end{align*}
\end{proof}

\begin{proposition}\label{goodboxprop}
Let $\alpha<1$.  There exist $M$ and $k$ with $M>2k$ such that for
all $p$ sufficiently close to 1, for every edge $(\bu,\bv)$,
\begin{equation}\label{wholeway}
\PP_p(\text{\rm edge $(\bu,\bv)$ is good})\geq \alpha.
\end{equation}
\end{proposition}

\begin{proof}
Pick $\beta\in(\alpha,1)$, and fix $M,k$ according to
Proposition \ref{goodboxprop2}.  Since the event that an edge is
good depends only on the initial states in a finite box, it is a
polynomial in $p$ and therefore continuous.  Thus the result
follows from Proposition \ref{goodboxprop2}.
\end{proof}

\begin{proofof}{Theorem \ref{maintheorem}}
Recall that the critical probability for oriented percolation
on $\ZZ^2$ is strictly less than 1 (see \cite{durrett} or \cite{g2}).
By the results of \cite{lss},
if $\alpha$ is sufficiently close to
1 then any 30-dependent bond percolation process on $\ZZ^2$
satisfying (\ref{wholeway}) stochastically dominates a Bernoulli
percolation process which is super-critical for oriented
percolation on $\ZZ^2$.

Therefore by Proposition~\ref{goodboxprop} and Lemma~\ref{mdeplemma}, we
may choose $M,k$ such that if $p$ is sufficiently close to 1, the event
that there is an infinite path of good renormalized edges starting from
$V_{(0,0)}$, oriented in the positive directions of both coordinates,
occurs with positive probability. On this event, there is an infinite
blocking path starting at $(0,0)$, so by Lemma \ref{blocking} we have
\[
  \PP_p(\text{there is a car which never moves at $(0,0)$})>0.
\]

Now consider any site $\bz$.  By translation invariance and
ergodicity, it follows from the above that almost surely there are
cars which never move at $\bz+(r,0)$ and $\bz+(0,s)$ for some (random)
$r,s\geq 0$.  This implies that any car initially at $\bz$ moves at
most $\max\{r,s\}$ times, while the state of $\bz$ changes at most
$2(r+s)$ times.
\end{proofof}

\begin{proofof}{Lemma~\ref{target}}
We start by giving an outline of the proof. Given a ``target''
$\by$, we will algorithmically construct a blocking path
$\bz^0,\bz^1,\ldots$ starting at $\bz^0=(0,0)$. If we use only
steps of types (i),(ii) in the definition of a blocking path, we
obtain a unique random path with asymptotic direction $(1,1)$. If
we also allow steps of types (iii),(iv), then at a positive
proportion of steps we have a choice of which direction to move.
By always choosing the direction which moves closer to the target
we are exponentially unlikely to miss the target by much, provided
that the target is within a cone determined by the typical slopes
that would result from choosing to go always up or always down.

We now present the details.  For simplicity, we will only allow
choices at alternate steps. Let $\bz^0=(0,0)$.  Suppose that a
blocking path $\bz^0,\ldots ,\bz^m$ has been constructed, and
suppose that $\bz^m$ lies on the diagonal line $z_1+z_2=2n$.  We
will extend the blocking path by one or two sites to some site on
the line $z_1+z_2=2n+2$.

Suppose first that
$\sigma(\bz^m)=\;\rta$, and consider the following cases:
\begin{itemize}
\item[(1)]
If $\sigma(\bz^m+(1,0))=\;\upa$ we set $\bz^{m+1}=\bz^m+(1,0)$ and
$\bz^{m+2}=\bz^m+(1,1)$.
\item[(2)]
If $\sigma(\bz^m+(1,0))=\sigma(\bz^m+(1,-1))=\;\rta$ we set
$\bz^{m+1}=\bz^m+(1,0)$ and $\bz^{m+2}=\bz^m+(2,0)$.
\item[(3)] If $\sigma(\bz^m+(1,0))=\;\rta$ and
$\sigma(\bz^m+(1,-1))=\;\upa$ we have a choice: we can set either
\begin{itemize}
\item[(a)]
$\bz^{m+1}=\bz^m+(1,0)$ and $\bz^{m+2}=\bz^m+(2,0)$
\item[or (b)]
$\bz^{m+1}=\bz^m+(1,1)$ (using a blocking path step of type (iii)).
\end{itemize}
We choose (a) if $z^m_1-z^m_2<y_1-y_2$, otherwise (b).
\end{itemize}

Thus we take the na\"{\i}ve path (using steps of types (i) and (ii)) unless
a step of type (iii) is possible and it moves us closer to $\by$ than the
alternative.

On the other hand if $\sigma(\bz^m)=\;\upa$ then $\bz^{m+1}$ (and possibly
$\bz^{m+2}$) are determined in an identical way, but interchanging the
roles of the two coordinates, and of $\upa,\rta$. In particular, in the
equivalent of case (3) above we extend the blocking path to
$\bz^{m+2}=\bz^m+(0,2)$ if $z^m_1-z^m_2>y_1-y_2$, and to
$\bz^{m+1}=\bz^m+(1,1)$ otherwise.

The above construction evidently yields a blocking path
$\bz^0,\bz^1,\ldots$. Suppose for the moment that $y_1+y_2$ is even. For
each $n$, let $\bz^{r(n)}$ be the site at which the blocking path
intersects the line $z_1+z_2=2n$, and let
$W_n=\big|(z^{r(n)}_1-z^{r(n)}_2)-(y_1-y_2)\big|/2$. It is straightforward
to check that $(W_n)_{n\geq 0}$ is a Markov chain with transition
probabilities
\[
\begin{array}{lll}
 P_{j,j-1}=1/4, & P_{j,j}=5/8, &  P_{j,j+1}=1/8
\quad\text{for $j\geq 1$}; \\
 P_{0,0}=3/4,  & P_{0,1}=1/4. & \\
\end{array}
\]
Thus $(W_n)_{n\geq 0}$ is a random walk on the natural numbers with drift
$-1/8$, and a reflecting boundary condition at $0$. To conclude we use the
following claim.

\begin{claim} \label{C:markov}
For the above Markov chain $(W_n)$, there exists $c_1>0$ such that for any
$N>9r$ and any $k$,
\[
  \PP(W_N>k \mid W_0=r) \leq e^{-c_1k}.
\]
\end{claim}

Assuming the claim we argue Lemma~\ref{target} as follows. If $y_1+y_2$ is
even, then the lemma follows from the claim immediately. If $y_1+y_2$ is
odd, then we apply the lemma first to $\by-(1,0)$ or $\by-(0,1)$ and $k-1$,
and note that that any finite blocking path may always be extended by one
site in direction $(1,0)$ or $(0,1)$.
\end{proofof}

\begin{proofof}{Claim~\ref{C:markov}}
Since the chain has increments at most 1, we have $W_N\leq r+N\leq
N/9+N<2N$.  Hence the probability in question is zero when $k>2N$, so
we may assume $k\leq 2N$.

Let $T$ be the first time $(W_n)$ hits $0$.  Before
 $T$, the increments are i.i.d. with mean $-1/8$, so
by the Chernoff bound we have $\PP(T>9r)\leq
e^{-c_2 N}\leq e^{-c_2 k/2}$.  Therefore, applying the strong
Markov property at $T$, the claim will follow if we can establish
for fixed $c_3>0$ and all $n\geq 0$ that $\PP(W_n>k\mid W_0=0)\leq
e^{-c_3k}$.  To check this, observe that we may couple $(W_n)$
with a stationary copy $(\widetilde{W}_n)$ in such a way that
$W_n\leq \widetilde{W}_n$ for all $n$, then note that the
stationary distribution has exponentially decaying tail.
\end{proofof}

\paragraph{Remarks.}
An alternative proof of Proposition \ref{goodboxprop2} involves
considering only blocking paths from the two endpoints of $V_\bu$
(rather than all $2k+1$ elements), and noting that blocking paths
cannot cross without intersecting.  (This argument does not extend to
higher dimensions).

Lemma \ref{target} in fact holds with the improved slope 3/2 (rather
than 9/8); this may be shown by allowing choices at all possible
steps rather than just alternate steps.

Experiments suggest that infinite blocking paths exist whenever $p>0.95$.

\section{Extensions}

At the core of the proof is a comparison of the collection of blocking
paths to super-critical oriented percolation. Since percolation is
relatively robust to variations in the model, it is not surprising
that our result holds for several other natural models. We present
several of these.

The proofs of the following theorems follow the same basic
argument as for Theorem~\ref{maintheorem}. Each of the variants
differs in some part of the proof, and so we only indicate the
changes that need to be made.  For simplicity we do not formulate
a model encompassing all extensions simultaneously.

\paragraph{The finite torus.}
We consider the model in which $\ZZ^2$ is replaced with the
rectangle
 $\{1,\ldots,m\}\times\{1,\ldots,n\}$ with periodic boundary
conditions.  Thus, a car moving East from $(m,i)$ re-appears at
$(1,i)$, while a car moving North from $(j,n)$ re-appears at
$(j,1)$.  Our aim is to prove Theorem \ref{torus}.

\begin{figure}
\centering \resizebox{!}{1.5in}{\includegraphics{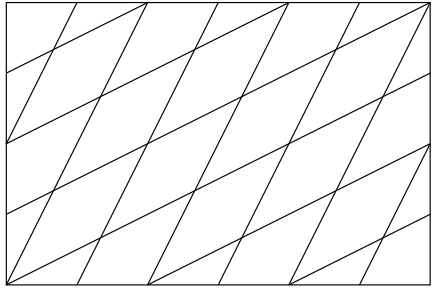}}
\caption{The renormalized lattice on a torus. Here it is the skew
torus ${\mathbb T}((6,-3),(-2,4))$.} \label{skew}
\end{figure}

We will use the following definition in constructing a renormalized lattice
on the torus. For linearly independent vectors $\ba,\bb\in\ZZ^2$, the {\bf
  skew torus} ${\mathbb T}(\ba,\bb)$ is the directed graph obtained from
the oriented square lattice by identifying vertices $\bx,\by\in\ZZ^2$
whenever $\bx-\by=s\ba+t\bb$ for some $s,t\in\ZZ$ (and identifying the
corresponding edges). See Figure~\ref{skew} for an illustration. The proof
of Theorem~\ref{torus} depends on the following lemma, which we prove by
standard percolation methods.

\begin{lemma}
\label{skew-perc} Let $q$ exceed the critical probability for
oriented bond percolation on $\ZZ^2$.  For any $\ba,\bb\in\ZZ^2$,
asymptotically almost surely as $r\to\infty$ with $r\in\ZZ$,
bond percolation with parameter $q$ on the skew torus ${\mathbb
T}(r\ba,r\bb)$ contains an open oriented cycle.
\end{lemma}

\begin{proof}
The event that bond percolation on the skew torus contains an open
oriented cycle is increasing, and it is quasi-symmetric; more
precisely, it is invariant under a group of permutations of the edges
of the skew torus having two transitivity classes (the horizontal
edges and the vertical edges).  By the Friedgut-Kalai sharp threshold
theorem (see Theorem 2.1 and the comment following Corollary 3.5 in
\cite{friedgut-kalai}), it therefore suffices to prove that for any
$q$ as described, the probability in question is bounded away from $0$
as $r\to\infty$.  (See \cite{bollobas-riordan} for another application of
\cite{friedgut-kalai} to percolation).

First consider oriented bond percolation with parameter $q$ on $\ZZ^2$.
We write $\bx\to\by$ for the event that there is an open oriented
path from $\bx$ to $\by$. We claim that
\begin{equation}
\label{diagonal}
\inf_{r\geq 1} \PP\big((0,0)\rightarrow(r,r)\big)>0.
\end{equation}
To check this, write $\theta=\PP((0,0)\to\infty)$ ($>0$), and note
that
\[
  \PP\big((0,0)\to(r,r)+(s,-s) \text{ for some } s\big) \ge \theta.
\]
Hence by symmetry,
\[
  \PP\big((0,0)\to(r,r)+(s,-s) \text{ for some }s \ge 0\big) \ge \theta/2,
\]
and similarly
\[
  \PP\big((s,-s)\to(r,r) \text{ for some }s \ge 0\big) \ge \theta/2.
\]
On the intersection of the last two events we have $(0,0)\to(r,r)$,
since the two directed paths must intersect. Therefore by the
Harris-FKG inequality (see \cite{harris} or \cite{g2}) we have
$\PP\big((0,0)\rightarrow(r,r)\big) \geq (\theta/2)^2$, establishing
(\ref{diagonal}).

Let $\ell$ be the smallest positive integer
such that $(0,0)=(\ell,\ell)$ in the skew torus ${\mathbb
T}(\ba,\bb)$ ($\ell$ is at most the number of vertices
in ${\mathbb T}(\ba,\bb)$). Now consider
bond percolation with parameter $q$ on ${\mathbb T}(r\ba,r\bb)$,
and let $A$ be
the event that
\[
  (0,0) \to (r,r) \to (2r,2r) \to \dots \to (\ell r,\ell r).
\]
By the Harris-FKG inequality
we have $\PP(A) \geq \gamma^\ell$, where $\gamma$ is the infimum in
(\ref{diagonal}).
And clearly on
$A$ there is an open oriented cycle.
\end{proof}

\begin{proofof}{Theorem \ref{torus}}
  First note that the existence of a cyclic blocking path including both
  horizontal and vertical steps is sufficient to ensure that no car moves
  infinitely often. The proof goes through as on $\ZZ^2$ except that we
  need to adjust the geometry of the renormalized lattice, which will have
  the structure of a skew torus. The process of good edges will still be
  30-dependent and the probability of an edge being good will be uniformly
  large provided the slopes of the renormalized edges lie strictly within a
  certain interval, and provided $M$ and $k$ are large enough; indeed $M$
  and $k$ may vary from edge to edge. Consider a sequence of tori of
  dimensions $m_k,n_k$ as in Theorem \ref{torus}. For $k$ sufficiently
  large we may construct a sequence of renormalized lattices subject to the
  above restrictions and with graph structure of ${\mathbb T}
  (r_k\ba,r_k\bb)$, where $r_k\to\infty$. The result then follows from 
  Lemma~\ref{skew-perc}, since an oriented cycle in the renormalized
  lattice yields the required cyclic blocking path.
\end{proofof}

\paragraph{Higher dimensions.}
Consider a variant model on $\ZZ^d$ in which each
non-empty site is occupied by a car facing in one of the $d$
directions. At times congruent to $i$ modulo $d$, all the cars
facing in direction $i$ advance if the place ahead of them is
empty. Many of the conjectures for the 2-dimensional model appear
reasonable in this case as well.  Define $\PP_p$ to be the
probability measure in which initially each site has a car with
direction $i$ with probability $p/d$, and is empty otherwise.

\begin{thm} \label{T:high_dim}
For the model on $\ZZ^d$ with any $d\ge2$, there exists some $p_1
= p_1(d) < 1$ such that for $p\geq p_1$, almost surely no car
moves infinitely often.
\end{thm}

\begin{proof}
The proof is very similar to that of the two-dimensional case.
Suppose more than one car is directly blocked by a car at $\bz$.
If the car at $\bz$ moves, than the order at which cars advance
dictates which of the blocked cars will enter $\bz$.  This allows
us to generalize the notion of a blocking path, and it is easy to
see that there is some fixed positive probability of being able to
continue a blocking path in any given direction.

The argument now continues as for $\ZZ^2$.  The probability of having
no path from $\bx$ to a neighbourhood of $\by$ inside a sufficiently
narrow cone is exponentially small, and the renormalization argument
applies.
\end{proof}

\paragraph{Biased initial conditions.}
Let $\PP_{\theta,p}$ be the probability measure on initial
configurations in which
$\PP_{\theta,p}\big(\sigma(\bz)=0\big)=1-p$ and
$\PP_{\theta,p}(\sigma(\bz)=\;\rta\!)=\theta p$ and
$\PP_{\theta,p}(\sigma(\bz)=\;\upa\!)$ $=(1-\theta)p$ for each
site $\bz$, and the states of different sites are independent.

\begin{thm} \label{T:bias}
  For any $\theta\in(0,1)$ there exists $p_1 = p_1(\theta) < 1$ such
  that for $p\geq p_1$ we have that $\PP_{\theta,p}$-a.s.\ no car moves
  infinitely often.
\end{thm}

\begin{proof}
  The proof of Theorem~\ref{maintheorem} adapts to this case as
  well. The maximum and minimum typical slopes of blocking paths are
  altered, and are not generally symmetric about the diagonal. A
  renormalized lattice spanned by two vectors inside the reachable
  cone can still be constructed.
\end{proof}

\paragraph{Random moves.}
Another interesting modification is to replace the deterministic
evolution of the model by a random mechanism.  In particular, suppose
each car attempts to move forward at the times of a Poisson process
of unit intensity, where different cars have independent Poisson
processes.

\begin{thm} \label{T:rand_time}
For the Poisson model there exists $p_1 < 1$ such that for any
$p\geq p_1$, almost surely no car moves infinitely often.
\end{thm}

\begin{proof}
  Consider a location corresponding to a step of type (iii) or (iv) in
  a blocking path. This involves a local configuration where two cars
  are directly blocked by a third car at some $\bz$. With
  deterministic evolution it is determined from the directions of the
  cars which of the two will advance to $\bz$ (thereby blocking the
  other), should $\bz$ become empty.  With the random moves this is
  not determined just by the directions. Clearly each of the two is
  equally likely to advance into $\bz$ before the other, independently
  of what happens at other locations where such a configuration
  exists.

  Thus we can toss an independent coin in advance at each such
  location, where the results of these coin tosses tell us which of
  the locations allow for a branching in the blocking paths and which
  do not. Since each potential branching point is retained with
  probability $1/2$ independently of all others, the blocking paths
  still form a super-critical process, and the proof goes through.
\end{proof}

\section*{Acknowledgements}
We thank Ehud Friedgut, Yuval Peres and Raissa D'Souza for
valuable conversations.  We thank the anonymous referee.

\bibliography{mybib}

\bigskip \noindent
{\bf Omer Angel}: \texttt{angel(at)math.ubc.ca}\\
{\bf Alexander E. Holroyd}: \texttt{holroyd(at)math.ubc.ca}\\
Department of Mathematics \\
University of British Columbia \\
Vancouver, BC V6T 1Z2, Canada

\bigskip \noindent
{\bf James B. Martin}: \texttt{James.Martin(at)liafa.jussieu.fr} \\
CNRS and Universit\'e Paris 7\\
2 place Jussieu (case 7014)\\
75251 PARIS Cedex 05, France

\end{document}